\RequirePackage{fix-cm}
\documentclass[smallextended]{svjour3}       
\smartqed  
\usepackage{graphicx}

\usepackage{amsmath}
\usepackage{amsfonts}
\usepackage{amssymb}
\newtheorem{Teo}{\textsc{Theorem}}
\newtheorem{Lemma}{{\textsc Lemma}}
\newtheorem{Prop}{\textsc{Proposition}}
\newtheorem{Cor}{\textsc{Corollary}}

\newcommand{\QED}{\hfill $\square$\bigskip}
\newcommand{\pf}{{\noindent\bf Proof.\ \ }}

\newcommand{\semidir}{\kern-2pt>\kern-6pt\triangleleft}
\newcommand{\N}{{\Bbb N}}
\newcommand{\Z}{{\Bbb Z}}
\newcommand{\Q}{{\Bbb Q}}

\newcommand{\J}{{\cal J}}
\newcommand{\p}{{\varphi}}
\newcommand{\e}{{\epsilon}}%

\newcommand{\TI}{{${\bf T}\kern-10pt\sim$}}

\newcommand{\iso}{\simeq}

\newcommand{\sn}{{\rm \kern1pt sn \kern1pt}}
\newcommand{\an}{{\rm \kern2pt an \kern1pt}}
\newcommand{\cf}{{\rm \kern2pt cf \kern1pt}}
\newcommand{\f}{{\rm \kern2pt f \kern1pt }}
\newcommand{\ff}{{\kern1pt \rm \bar{f}\kern1pt }}
\newcommand{\C}{\rm \kern2pt c \kern1pt}
\newcommand{\cn}{\rm \kern2pt cn\kern1pt}
\newcommand{\ssim}{\rm \kern2pt \approx}

\usepackage{graphicx}

\usepackage{amsmath}
\usepackage{amsfonts}
\usepackage{amssymb}
\journalname{ }
\begin{document}



\title{
On the ring of  inertial endomorphisms of an abelian group}

\author{Ulderico Dardano\  \ -\
Silvana Rinauro }

\institute{Ulderico Dardano \at
             Dipartimento di Matematica e
Applicazioni ``R.Caccioppoli'', Universit\`a di Napoli ``Federico~II'', Via Cintia - Monte S. Angelo, I-80126 Napoli, Italy. \email{dardano@unina.it}
           \and
           Silvana Rinauro \at
Dipartimento di Matematica, Informatica ed Economia, Universit\`a della
Basilicata, Via dell'Ateneo Lucano 10,
I-85100 Potenza, Italy. \email{silvana.rinauro@unibas.it}}

\date{Received: date / Accepted: date}

\maketitle


\maketitle

\vskip-0.7cm
\noindent{\textbf{Abstract.}\ \ \ An endomorphisms $\p$ of an abelian group $A$ is said inertial if each subgroup $H$ of $A$ has finite index in $H+\varphi (H)$. We
study the ring of inertial endomorphisms of an abelian group. We obtain a satisfactory description modulo
the ideal of finitary endomorphisms.
Also the corresponding problem for vector spaces is considered.

\noindent{{\bf 2000 Mathematics Subject Classification}:\ 20K30 (16S50).



\medskip
\section{Introduction and linear version}

If $\p$ is an endomorphism of an abelian group $A$ and $H$ a subgroup of $A$, then in \cite{DGBSV} $H$ is called a \emph{$\p$-inert subgroup} iff $H$ has finite index in $H+\p(H)$. Obviously, finite, finite index and $\p$-invariant subgroups are $\p$-inert.
  Passing to a \lq\lq global condition", we have -on one side- the
notion of fully inert subgroup as in \cite{DSZ}, that is $\p$-inert w.r.t. any endomorphism $\p$. Motivation for studying {fully inert}
subgroups comes from the investigation of the dynamical properties of
an endomorphism of an abelian group (see \cite{DGSZ}, \cite{DGBSV}).
{Fully inert}
subgroups of an abelian group $A$ have been recently studied in \cite{DGBSV2}, \cite{DSZ},   \cite{GSZ}, \cite{GSZ2}  in cases when $A$ is a divisible group, a free group, a torsion-free module over the ring of p-adic integers, a $p$-group, resp.

On the other side, we call \emph{inertial an endomorphism with respect to which each subgroup is inert}. This property has been studied in \cite{DR1} and \cite{DR2} in connection with the study of inert subgroups of groups (see \cite{BKS}, \cite{DGMT}, \cite{R}). Recall that in non-abelian group context,  a subgroup is said to be \emph{inert} if it is commensurable with its conjugates (that is w.r.t. inner automorphisms), where two subgroups $H$ and $H_1$ of a group are told commensurable iff $H\cap H_1$ has finite index in both $H$ and $H_1$.

 Clearly the definition of inertial endomorphism of an abelian group can be regarded as a local  condition generalizing both notions of {power} and {finitary} endomorphism. Recall that an endomorphism is \emph{power} if it maps each subgroup to itself and \emph{finitary} if it has finite image.

In \cite{DR2} we have proved the following Fact (and given a fruithful description of inertial endomorphisms as reported in the Characterization Theorem below).

\medskip
\noindent
{\bf Fact} \emph{Inertial endomorphisms of any abelian group $A$ form a subring
 $IE(A)$ of the full ring $E(A)$ of endomorphisms of $A$. It contains
 the ideal $F(A)$ of finitary endomorphisms and
 many "multiplications" (see below). Moreover inertial endomorphisms of an abelian group commute modulo finitary ones, that is $IE(A)/F(A)$ is a commutative ring}.

\medskip

In this paper \emph{we study the ring $IE(A)$ of inertial endomorphisms}, while the group of invertible inertial endomorphisms
has been studied in \cite{DR3} along the same lines as in \cite{W} and \cite{BS1}.

We first note that the corresponding problem for vector spaces has an easy solution.

\begin{Teo}\label{TEOREMA_LINEARE--} Let $\p$ be an endomorphism of an infinite dimension  vector space $V$. Then each subspace $H$ has finite codimension in $H+\p(H)$  iff $\p$ acts as a scalar multiplication on a finite codimension subspace of $V$.

On the other hand, for each subspace $H$ the codimension of $H \cap \p(H)$ in $H$ is finite  iff
$\p$ acts as multiplication by a non-zero scalar on a finite codimension subspace.
\end{Teo}

 \begin{Cor} The ring of endomorphisms $\p$ of a vector space $V$ such that each subspace $H$ has finite codimension in $H+\p(H)$ is the sum of the ring of scalar multiplications and the ring of finitary endomorphisms.
  \end{Cor}

These statements apply also to the above defined ring $IE(A)$ in the case $A$ is elementary abelian, clearly. A corresponding statement holds for \lq\lq elementary" inertial endomorphisms (see section \ref{FM}). To treat the periodic case, we introduce \emph{mini-multiplications}, see sections \ref{mini-moltiplicazioni} and \ref{periodico}.
 In the general case, the picture can be rather complicated, see our Main Result Theorem \ref{TeorIENDmisto}. Further relevant instances of (uniform) inertial endomorphisms  are introduced in section \ref{UI} by Proposition \ref{RI} and Corollary \ref{G}.

\medskip
\noindent \textbf{Proof of Theorem \ref{TEOREMA_LINEARE--}}
We use an argument similar to one we used in the proof of Proposition 5 of \cite{DR2}.
Fix $\p$ and suppose that either for each $H\le A$ the group $H+\p(H)/H$ is finite or
 for each $H\le A$ the group $H/H\cap \p(H)$ is finite.
By contradiction, assume $\p$ is scalar multiplication on no quotient by a finite dim subspace.
We first prove that \emph{if $W$ is a finite dim subspace, then $W$ is contained in a finite dim $\p$-invariant subspace}. Clearly we can assume that $W=Ka$
has dimension $1$. Consider then the homomorphism
$$\Phi:K[x]\mapsto V$$ mapping $1\mapsto a$ and $x\mapsto\p(a)$.
If $\Phi$ is injective,
we can replace $im\Phi$ by $K[x]$ and $\p$ by multiplication by $x$. Then both $H:=K[x^2]$ and $\p(H)=xH$ have infinite dimension,
while $H\cap xH=0$, a contradiction. Therefore $\Phi$ is not injective and
$im\Phi$ is the wished subspace.

Now we can prove that:
\emph{for all finite dimension subspaces $W\le V$
 such that $W\cap \p(W)=0$ there exists a subspace $W'> W$ with finite dimension such that $$W'\cap \p(W')=0\ \ {\rm and } \ \ \p(W')>\p(W).$$}To see this note that if $Z\ge W$ is a finite dim $\p$-invariant subspace, as $\p$ does not act as a scalar multiplication on $A/Z$,
 we can choose $a\in
V$ such that
 $\p(a)\not\in Ka+Z$ and define $W':=Ka+W$.
If now
 $b\in W'\cap \p(W')$,  then $\exists h,k\in K$, $\exists c,c_0\in W$ such
 that $b=ha+c=k\p(a)+\p(c_0)$. Thus $k\p(a)\in Ka+
Z$ while $\p(a)\not\in Ka+Z$. Therefore $k=h=0$. It follows $b=c=\p(c_0)\in W\cap
 \p(W)=0$, as claimed.
\smallskip

 Finally, starting at $W_0=0$, by recursion we
 define $W_{i+1}:=W_i'$ and $W_\omega:=\cup_i W_i$. We get that both $W_\omega$ and $\p(W_\omega)$ have
 infinite dimension and $W_\omega\cap \p(W_\omega)=0$, the wished contradiction. The statement follows now easily.
\qed

\section{Notation and statement of Main Result}
To state our main result,
we need to consider some relevant invariants of the group and to introduce some \emph{ad hoc definitions}. Note that in this paper $A$ always stands for an abelian group (in additive notation).
For undefined notation on abelian groups we refer to \cite{F}. In particular
we denote by $Q_p^*$  the ring of $p$-adic integers and by $\J$ the ring $\prod_p Q_p^*$ which is meant to act componentwise on any periodic abelian group. For $m,n\in\Z$, if we consider the fraction $m/n\ \in\Q$, we always mean $m$ and $n\ne0$ are coprime. We write  $\pi(n)$ for the set of prime divisors of $n\in\N$ and
$\pi'$ for the complement of a set $\pi$ of primes.
Denote by $A_\pi$ the $\pi$-component of $A$ and $A[n]:=\{a\in A\ |\ na=0\ \}$. As usual, if $n\in\N$ annihilates $A$, that is $nA=0$, we say that $A$ is \emph{bounded} by $n$. Further, we say that $A$ is bounded, if it is bounded by some $n$ and the least such $n$ is called \emph{the bound} of $A$. If $pA=A$ (resp. $A$ is periodic and $A[p]=0$) for each $p\in\pi$ say that $A$ is $\pi$-divisible (resp. a $\pi'$-group).

We recall that in \cite{DR2}, we called \emph{multiplications} of an abelian group $A$ either the actions on $A$ of a subring of $\Q$ or, when $A$ is periodic, the above action of $\J$.
\emph{Multiplications form a ring $M(A)$}. If $A$ is periodic, clearly $M(A)\iso\ \prod_p\ Q_p^*/(p^{e_p})$, where the product is taken on all primes and $p^{e_p}$ is the bound of the $p$-component $A_p$ of $A$ or $0$ if  $A_p$ is unbounded. If $A$ is non-periodic, then $M(A)\iso\Q^{\pi}$ where $\pi$ is the largest set of primes such that $A$ is a $\Q^\pi$-module, that is $A$ is $p$-divisible with no elements of order $p$, $\forall p\in\pi$. Here, as usal, $\Q^{\pi}$ is the ring
of rationals whose denominator is a $\pi$-number, that is divided by primes in $\pi$ only.

 Recall also that if $A$ is any abelian group there is a free abelian subgroup $F$ such that  $A/F$ is periodic. The rank of $F$ coincides with the torsion-free rank $r_0(A)$, that is the rank
of the torsion free group $A/T$, where $T=T(A)$  denotes the torsion subgroup of $A$, as usual. In Proposition 1 of \cite{DR2} we noticed that \emph{when $A$ is not periodic multiplications which are not by an integer are inertial iff the underlying abelian  group $A$ has finite torsion-free rank}. For abelian groups with infinite torsion-free rank, case $(a)$ in Characterization Theorem below. Thus
\emph{we will be mainly concerned with  groups with finite torsion-free rank, say FTFR}.

When $A$ has FTFR all subgroups $F$ as above are commensurable. Fix one of them and
define the \emph{sequence of essential $p$-bounds} $(\e_p)$
such that $\e_p$ is the min $i$ such that $\sum_{j>i}f_j<\infty $ where $f_j$ is the \emph{Ulm-Kaplanski invariant} of the $p$-component $(A/F)_p$
and $\e_p=\infty$ if the $p$-component of $A/F$ is unbounded. Clearly $p$ runs on all primes.
In other words, $\e_p=\infty$ or $\e_p$ is the smallest $\e$ such that $p^\e(A/F)_p$ is finite.
Clearly the sequence $(\e_p)$ is independent of $F$.

We also consider the \emph{sequence of $p$-bounds} $(e_{p})$, where $p^{e_p}$ is
the bound of $(A/F)_{p}$, when this is bounded, or $p^{e_p}=\infty$, otherwise. This sequence
depends on the choice of $F$. Clearly $\e_{p}\le{e_p}$ for each prime $p$. However sequences $(e_{p})$'s corresponding to different $F$ are definitely equal and coincide on entries which are infinite. Then we may consider an equivalence relation such that
the class of the sequence  $(e_{p})$ depends on $A$ only.

Thus we define (independently of $F$ and of the appearing $e_p$'s) the ring
$${\cal H}(A):=\dfrac{\prod_p\ Q_p^*/(p^{e_p})  }{\bigoplus_p\ (p^{\e_p})/(p^{e_p})}
\iso \dfrac{\prod_p\ Q_p^*}{\bigoplus_p\ (p^{\e_p})\ +\ \prod_p\ (p^{e^p})}  $$
where $(e_p)$ and $(\e_p)$ are sequences of $p$-bounds and essential $p$-bounds resp, $(p^x)$ is ideal generated by $p^x$ in  the ring $Q_p^*$, where $p^\infty:=0$. Also denote by:\\
- $\pi_*(A)$ the set of primes $p$ such that $A_p$ is bounded and $A/A_p$ is $p$-divisible.\\
-  $\pi_c(A)$ the set  of primes for which for some $F$ as above the $p$-component of $A/F$ is \emph{critical}, that is with shape $B\oplus D$ with $B$ infinite but bounded
and $D\ne 0$ divisible with finite positive rank. Clearly $\pi_c(A)$ is independent of $F$, as $A$ has FTFR.
\medskip

As we shall often face some particular decomposition of an abelian group let us introduce further  notation. If $A=B\oplus C$ and there is a \emph{finite} set $\pi$ of primes, such that $B$ is  a \emph{bounded} $\pi$-subgroup and the $\pi$-component $C_{\pi}$ of $C$ is \emph{divisible with finite rank} we write\\ \centerline{$A=B \oplus_\pi C$}\\
and call \emph{mini-multiplications} endomorphisms of shape $n\oplus 0$ (i.e. acting as $n\in\Z$ on $B$ and annihilating $C$). For details see section \ref{mini-moltiplicazioni})

Further, it is easy to verify that
\emph{if $\pi\subseteq \pi_*(A)$,
then $A=B \oplus_\pi C$  for $B=A_\pi$  and subgroup $C$ is a $\Q^{\pi}$-module which is uniquely determined and fully invariant}, clearly.
In such a condition, we say that $\p:=n\oplus_\pi \alpha$ is a \emph{semi-multiplication} where $\p$ is the  multiplication by $n\in\Z$ on $B$ and by
$\alpha$ on $C$, where either $\alpha\in \J$ or  $\alpha\in \Q^{\pi}$ according to
$C$ is periodic or not. Denote by $SM(A)$ the ring formed by the semi-multiplications. It is  \emph{is a central subring of $E(A)$}.

Finally, say that an inertial endomorphism $\p$ is \emph{uniform} iff there is a $\p$-series $F\le A_0\le A$ such that $A/A_0$ is finite,  $\p$ is multiplication on the periodic group $A_0/F$ and  $\p=0$  on the free abelian subgroup $F$. Denote by $UI(A)$ the \emph{subring of uniform inertial endomorphisms} (see section \ref{UI}).

Now we can describe the ring $IE(A)$ of inertial endomorphisms in terms of above invariants and others that we will introduce below. Symbols of (direct) sum are to be read in the additive group of  $IE(A)$.

\begin{Teo}\label{TeorIENDmisto}\emph{ \bf (Main Theorem)}\
Let $A$ be an abelian group. \\
a) If $A$ has not FTFR,
then $IE(A)=M(A)\oplus F(A)$, where $M(A)\iso\Z$ is the ring of multiplications by integers.\\
\smallskip
b) If $A$ has FTFR, then
$$IE(A)=SM(A)+UI(A)+N$$
where:\\
- $SM(A)$ is the (central) subring of \emph{semi-multiplications} of $A$;\\
-  $UI(A)$ is the subring of \emph{uniform inertial} endomorphisms of $A$;\\
-  $N\iso \oplus_{p\in\pi_c(A)} \Z(p^{c_p})$ is a subring of inertial \emph{mini-multiplications} of $A$.\\
\smallskip
Further, we have:\\ \smallskip
$1)$  $UI(A)+N$ is the ideal of inertial endomorphisms of $A$ with periodic image;\\
$2)$  $\dfrac{IE(A)} {UI(A)+N}\iso \Q^{\pi_*(A)}$ as rings (provided $A$ is non-periodic);\\
$3)$  $\dfrac{UI(A)}{F(A)}$ is isomorphic to a subring of  ${\cal H}(A)$\ \ .
\smallskip\\
c) In particular, if $A$ is periodic, then:\\
$i)$\ \ \   $IE(A)=M(A)+F(A)+N;$\\
$ii)$ \ \ $SM(A)=M(A)$ \ \ and\ \ $UI(A)=M(A)+F(A)$;\\
$iii)$ $\dfrac{UI(A)}{F(A)}\iso{\cal H}(A)$ as rings.
\end{Teo}

For details on above $N$ and $c_p$'s see Lemma \ref{mini-moltiplicazioni}. By Corollary \ref{G} we will see that $UI(A)$ may be rather large even if $A$ is countable. Now we recall a characterization of inertial endomorphisms (see Theorem A and Proposition 5 in \cite{DR2}).

\bigskip
{\small

\noindent\textbf{Characterization Theorem} (\cite{DR2}). {\it Let $\p_1,\dots,\p_l$ finitely many endomorphisms of an abelian  group $A$. Then each
$\p_i$ is inertial if and only if there is a finite index subgroup $A_{0}$ of $A$ such that one of (a) or (b) holds:
\noindent
$(a)$ each $\p_i$ acts as multiplication on $A_{0}$ by $m_i\in\Z$;

\noindent $(b)$ $A_{0}=B\oplus D\oplus C$ and exist finite sets of primes $\pi\subseteq \pi_1$ such that: \\
\phantom{m}  i)\ $B\oplus D$ is the $\pi_1$-component of $A_{0}$ where $B$ is bounded, $D$ divisible $\pi'$-group with finite rank.\\
 \phantom{m}
ii)\ $C$ is a $\Q^{\pi}[\p_1,\dots,\p_l]$-module, with a submodule $V\iso \Q^{\pi}\oplus\dots\oplus\Q^{\pi}$ (finitely many times) such that $C/V$ is a $\pi_1$-divisible $\pi'$-group.\\
 \phantom{m}
iii)\ each $\p_i$ acts by (possibly different) multiplications on $B$, $D$, $V$, $C/V$ where $\p_i$ is represented by $m_i/n_i\in\Q$ on $V$ and on all $p$-components of $D$ such that the $p$-component of $C/V$ is infinite and $\pi=\pi(n_1\cdots n_l)$.

Moreover, if $A$ is periodic, then $\p_1,\dots,\p_l$ are inertial iff\\
(FS) there is $m\in\N$ such that for each $X\le A$ there are subgroups $X_*,X^*$ which are $\p_i$-invariant ($\forall i$) and such that $X_*\le X\le X^*\le A$ and $|X_*/X^*|\le m.$
}

\medskip
Clearly if $A$ is torsion-free inertial then endomorphisms are multiplication by rationals.
When $A$ is periodic we have $V=0$ and in particular:\\
- if $A$ divisible and periodic then inertial endomorphisms are multiplication,\\
- if $A$ is reduced and periodic inertial endomorphisms are multiplication on a finite index subgroup of $A$.
}


\section{The ring $FM(A)$ of multiplication-by-finite endomorphisms}\label{FM}

In this section, by Theorem \ref{PF} we study $FM(A)$, a relevant subring of $IE(A)$, which might have interest in itself, as well. The following is easy to check.

\medskip \noindent
 \textbf{Fact} \emph{If  $\p\in E(A)$ acts as an inertial endomorphism on a finite index subgroup of $A$ then $\p$ is inertial on the whole $A$ indeed. The same happens arguing modulo a finite $\p$-invariant subgroup.}

 \medskip \noindent
We say that two endomorphisms are \emph{close} iff their difference is finitary.
We generalize Proposition 4 of \cite{DR2} to non-periodic groups and give a picture of the ring of endomorphisms which are close to multiplications. We give a definition.

Denote by $\pi_0(A)$ the set of primes $p$ such that $A_p$ is finite and $A/A_p$ is $p$-divisible. Clearly $\pi_0(A)\subseteq \pi_*(A)$.  If $\pi$ is a finite subset of $\pi_0(A)$ then $A=A_\pi\oplus_{\pi} C$. Notice that summands are fully invariant and uniquely determined. Call \emph{quasi-multiplication} of a non-periodic group $A$ those endomorphisms with shape $r\oplus_\pi m/n$ (with $r\in\Z$ and $m/n\in\Q^\pi$). Clearly these form a subring $QM(A)$ of $E(A)$. If $A$ is periodic set $QM(A):=M(A)$. In any case,  $QM(A)\subseteq SM(A)$ and $QM(A)$ is also in the center of $E(A)$.

 \begin{Teo}\label{PF} For an endomorphism $\p$ of an abelian group $A$ the following are equivalent:\\
MF)\ \ $\p$ acts by means of a multiplication on a finite index subgroup $A_0$ of $A$,\\
FM)\ $\p$ acts by means of a multiplication modulo a finite subgroup $A_1$ of $A$.

Moreover, endomorphisms with the above properties form a subring
$$FM(A)=F(A)+QM(A)$$
which is contained in $IE(A)$, provided $A$ has FTFR. Moreover,\\
$i)$ \ if $A$ is non-periodic, then
 $FM(A)/F(A)$ is naturally isomorphic to the ring $\Q^{\pi_0(A)}$ and
$F(A)\cap QM(A)$ consists of maps of type $r\oplus_\pi 0$ for a finite subset $\pi$ of  $\pi_0(A)$,\\
$i)$  if $A$ is a $p$-group, then $QM(A)=M(A)$; if $A$ is unbounded it holds
$F(A)\cap M(A)=0$; otherwise, if $e<\infty$ and $\e$ are the $p$-bound and the essential $p$-bound of $A$, resp., there is a natural ring isomorphisms $F(A) \cap M(A) \iso p^\e\Z/p^e\Z.$
\end{Teo}

\noindent Note that an endomorphisms of a periodic abelian group $A$ is (FM) iff it acts this way on finitely many components and by multiplications on all remaining ones, clearly. Note also that on $A=\Z(p)\oplus_p\Q^{\{p\}}$ we have $0\oplus_p {1}/{p}\in FM(A)\setminus(F(A)+M(A))$.

\medskip

\noindent \pf$(MF)\Rightarrow (FM)$ If $A$ is periodic, the statement is clear as there is  $\alpha\in\J$ such that $\p=\alpha$ on $A_0$ and we can consider $A_1:=im(\p-\alpha)$. Otherwise, $\p=m/n\in\Q$ on $A_0$, a $\Q^{\pi}$-module with $\pi:=\pi(n)$.  Thus $A_\pi$ is finite and we may assume $A_\pi=0$. Then note that $A/T$ is torsion free with a finite index $\pi$-divisible subgroup $A_0+T/T$. Then $A/T$ is $\pi$-divisible. Therefore $A$ is $\Q^{\pi}$-module. Again, $A_1:=im(\p-{m/ n})$ is the wished subgroup.

$(FM)\Rightarrow (MF)$  If $A$ is non-periodic and $\p=m/n$ on the $\Q^\pi$-module $A/A_1$ (again $\pi:=\pi(n)$), then $A_\pi\le A_1$ is finite and there a $\Q^\pi$-module $A_0'\le A$ such that $A=A_\pi\oplus A_0'$ and one can consider the endomorphism $\p_0:=0\oplus  m/n\in E(A)$ such that $im(\p-\p_0)\le A_1$. Thus  $A_0:=A_0'\cap ker(\p-\p_0)$ has finite index in $A$. On the other hand since $A_0'$ is $\pi$-divisible and $A_0'/A_0$ is finite we have
$A_0'/A_0$ is coprime to $n$ and so $A_0$ is $n$-divisible, thus $\p_0$-invariant.
Therefore $A_0$  is $\p$-invariant too and is the wished subgroup. The periodic case is clear, with $\p=\alpha\in\J$ on $A/A_1$ and $\p_{0}=\alpha\in M(A)$.

Moreover, let $\p_i$ act by multiplication on $A/A_i$ ($i=1,2$) with $A_3:=A_1+A_2$ finite.
As any subgroup of $T(A)/A_i$ is $\p_i$-invariant ($i=1,2$), such is $A_3$. Then $\p_1-\p_2$ and $\p_1\p_2$ act as multiplications on $A/A_3$.

Above arguments also show that $FM(A)=F(A)+QM(A)$ as $\p_0\in QM(A)$.
 From the Characterization Theorem above
it follows that $FM(A)\subseteq IE(A)$ when $A$ has FTFR.

Finally, when $A$ is non-periodic, the map $\p=r\oplus m/n\in FM(A)\mapsto m/n\in \Q^{\pi_0(A)}$ is the wished isomorphism and if
 $r\oplus_\pi m/n\in QM(A)$ is finitary, then it is $0$ on $A/T(A)$. Hence $m=0$. On the other hand, when $A$ is a $p$-group,
if $0\ne\alpha\in M(A) \cap F(A)$ we have that there
exists $i$ such that $ker\ \alpha=A[p^i]$ (clearly $p^i$ is the maximal power of $p$ dividing $\alpha$).
If $A[p^i]$ has finite index in $A$, then $A$ is bounded and $\e\le i$.
Conversely, if $p^\e$ divides $\alpha$ it is plain that $\alpha\in F(A)$.
\qed


\section {Mini-multiplications of an abelian group}

From the introduction section recall the following:

\noindent \textbf{Definition} A \emph{mini-multiplication} of an abelian group $A=B\oplus_{\pi}C$ is an endomorphism of shape $n\oplus 0$ i.e. acting as multiplication by $n\in\Z$ on the $\pi$-bounded subgroup $B$ and annihilating $C$, where the $\pi$-component of $C$ is divisible with finite total rank.

 By next statement we consider the isomorphism type $NM(A)$ of the
subring $N$ which appears in the statement of Theorem \ref{TeorIENDmisto}. It is the type of a ring of multiplications with bounded support of a fully-invariant distinguished periodic section of $A$ as well.

\begin{Lemma}\label{mini-moltiplicazioni} \  Let $A$ be an abelian group with FTFR and
$p_i^{c_i}$ be the bound of $A_{p_{i}}/div(A_{p_{i}})$ for each $i$ where $\pi_c(A)=:\{p_1,\ldots,p_i,\ldots\}$ is the set of critical primes of $A$.
Let $\pi_{i}:=\{p_1,\ldots,p_i\}$.

There are sequences $(B_{i})$ and $(C_{i})$ of subgroups such that $B_{i}$ is a bounded $p_i$-group $\forall i$ and\\ $(*)\ \ \ \ \ \ A=(B_{1}\oplus...\oplus B_{i}) \oplus_{\pi_{i}}C_{i}\ \ \ {\rm and}\ \ \  C_i=B_{i+1}\oplus C_{i+1}.$\\

Fixed above sequences, mini-multiplications $n\oplus_{\pi_{i}}0$ acting as $n\in\Z$ on $B_{1}\oplus...\oplus B_{i}$ and $0$ on $C_i$
(for some $i$) form a ring $N$ of inertial endomorphisms isomorphic to
$$NM(A):=\oplus_{p_i\in\pi_c(A)} \Z(p_i^{c_i}).$$
\end{Lemma}

\pf Note that the $p_{i}$-component of $A$ is the sum of a bounded subgroup and a divisible one. Therefore $A$ splits on it.
We define inductively wished sequences $(B_{i})$ and $(C_{i})$ by (choosing) $A=B_1\oplus C_1$ and $C_i=B_{i+1}\oplus C_{i+1}$. Then define  mini-multiplications  as in the statement.
They are inertial by the Characterization Theorem. Notice that unfortunately
they depend
on the choice of the two sequences.
Then for each $i$ and coset class $n_i\in \Z(p_i^{c_i})$
consider the mini-multiplication which is multiplication by $n_i$ on $B_i$ and $0$  on
other summands. This gives the wished isomorphism.
\qed

Let us highlight the role of mini-multiplications. Say that an endomorphism $\psi$ of an abelian group $A$ is \emph{bounded} iff $\psi(A)$ is bounded.

\begin{Lemma}\label{Lemma-moltiplicazioni-bounded} Let $\mu$ be a bounded multiplication of an abelian $A$. If $A$ is non-periodic, then $\mu=0$. If $A$ is a $p$-group $A$, then either  $\mu=0$ or $A$ is bounded.
\end{Lemma}

\pf If $A$ is non-periodic, consider the action of $\mu$ on $A/T(A)$. Otherwise assume $\mu\not=0$ and write $\mu=p^r\alpha$ where $\alpha$ is an invertible $p$-adic and $r\in\N_0$ and consider that  $p^rA=\mu(A)$, whence $A$ is bounded. \QED

\begin{Lemma}\label{Lemma-minimoltiplicazioni}
Let $\p$ be an endomorphism of $A$. If $\p$ acts as a mini-multiplication of type $m\oplus _{\pi}0$ on a subgroup of finite index  $A_0$ of $A$, then $\p$ is close to a mini-multiplication of type $m\oplus_{\pi} 0$ on the whole of $A$.
\end{Lemma}
\pf
Let  $\p=m\oplus_{\pi} 0$ on $A_{0}=B_{0}\oplus_{\pi}C_{0}$. Then $A_{\pi}=B\oplus D$, where $D\le C_{0}$ is divisible of finite rank and $B\ge B_{0}$ is bounded. Hence  $A=A_{\pi}\oplus C_{1}$ and $C:=D\oplus C_{1}$  is commensurable to $C_{0}$. It is now clear that $\p$ is close to the mini-multiplication $m\oplus_{\pi}0$ on $A=B\oplus_{\pi}C$.
\qed

\begin{Prop} Let $A$ be an abelian group.\\
a) If $A$ has not FTFR, bounded inertial endomorphisms are finitary. \\
b)  If $A$ has FTFR, any bounded inertial endomorphism is the sum of a mini-multiplication and a finitary endomorphism.
\end{Prop}

\pf Let $\p$ be a bounded inertial endomorphisms of $A$.
If $A$ has not FTFR, then $\p$ is multiplication on a finite index subgroup of $A$ and the statement follows as in Lemma \ref{Lemma-moltiplicazioni-bounded}.

If $A$ has FTFR then we are in case $(b)$ of the Characterization Theorem above. We use the same notation. Then $\p=m\in\Z$ on $B$ and $\p=0$ on $D$, as this is divisible. Let $\pi_2$ be the set of primes $p$ in $\pi(\p (C))$ (which is finite) such that the $p$-component of $C$ is bounded. Then $C=C_{\pi_2}\oplus C_1$. Also $\p=0$ on $V$. Moreover $C_1\cap V$ is $\p$-invariant  and $C/(C_1\cap V)$ is periodic. Then by (FS), there is a finite index $\p$-invariant  subgroup $C_0$ of $C_1$.

Note that if $p\in\pi(\p(C_0))$ then $p\not\in\pi_2$, hence $\p=0$ on
 the (unbounded) $p$-component of $C_0$, by Lemma \ref{Lemma-moltiplicazioni-bounded}.
 Therefore $\p=0$ on the whole $T(C_0)$. Then
  $\p(C_0)$ is an image of $C_0/T( C_0)$ which has finite rank. As $\p(C_0)$ is bounded, it is even finite. Thus $\p=0$ on a finite index subgroup $C'$ of $C_0$. Note that $C'$
is a $\pi_2$'-group and $C'_{\pi_1}$ is divisible with finite rank as $C_{\pi_1}$ is.  Thus
$\p$ is mini-multiplication on $(B+C_{\pi_2})\oplus_{\pi_1\cup \pi_2} C'$, which has finite index in $A$. By Lemma \ref{Lemma-minimoltiplicazioni}, $\p$ is close to a mini-multiplication. \QED


\section{Periodic case}\label{periodico}


The ring $IE(A)$ when $A$ is a  $p$-group is described by the following result,
which follows from the Characterization Theorem above. For the general periodic case note
that an \emph{endomorphism
of an abelian torsion group is inertial iff it is such on all primary
components and multiplication on all but finitely many of them.}

Recall that we say that an abelian $p$-group is \emph{critical} when $\pi_c(A)=\{p\}$, that is when the maximum divisible subgroup $D:=D(A)$ has positive finite rank and $A/D$ is bounded but infinite.

\begin{Prop}\label{ENDO_mod_F}
Let $A$ be an abelian $p$-group.\\
$i)$ \ If $A$ is non-critical then
$IE(A)= FM(A)$.\\
$ii)$ If $A$ is critical and $D:=D(A)$ then
$$IE(A)= FM(A) + N= M(A)+F(A)+N$$
where $N \iso M(A/D)\iso\Z(p^e)$ is a subring
of mini-multiplications and $p^{e}$ is the bound of $A/D$.

Moreover\ \ $M(A) \cap (F(A) + N)=0$\ and \ $FM(A)\cap N=F(A)\cap N=p^\e N$ where
 $\e$ is the essential $p$-bound of $A/D$.
\end{Prop}

\pf Let $\p\in IE(A)$. If $A$ is non-critical, $\p\in FM(A)$ by the Characterization Theorem.
If $A$ is critical, fix a decomposition $A=B\oplus_p D$ with $B$ bounded and $D=D(A)$ with finite rank. By  the same theorem there is a $\p$-invariant finite index subgroup $B_0$ of $B$ such that $\p$ acts as multiplication by some $n\in\Z$ on a finite index subgroup $C$ of $B_0$.

For each $p$-adic $\beta\in Q_p^*$ denote by $\bar \beta$ the mini-multiplication $\bar \beta=\beta\oplus 0$ on $A=B\oplus_p D$. Then if $\alpha\in Q_p^*$ represents the action
of $\p$ on $D$ (which is multiplication) we have $C+D\subseteq ker(\p-\alpha-\overline{(n-\alpha)})$ and so
$IE(A)= M(A) + F(A) + N$, where $N=\{\bar n\ |\ n\in\Z\}\iso \Z(p^e)$.

Further,
if $\alpha=\p_0+\bar n\in M(A) \cap (F(A) + N)$ (where $\p_0\in F(A)$ and $\bar n\in N$), then $\alpha=0$ on $D$. It follows $\alpha=0$ on the whole of $A$. By a similar argument $FM(A)\cap N=F(A)\cap N$. Finally, if the mini-multication $m\oplus 0$ on the decomposition $A=B\oplus D$ is finitary, then $m=0$ on $B[p^\e]$ and $p^\e$ divides $m$.
\QED

Recall that the description of
group of units of $IE(A)$ given in \cite{DR3} is rather complete when $A$ is periodic.




\section{The ring UI(A) of uniform inertial endomorphisms}\label{UI}

We consider now a further relevant ring of inertial homomorphisms,
which contains $F(A)$. \\
\textbf{Definition} An
inertial endomorphsm $\p$ is \emph{uniform and represented by $\alpha \in \J$ on $A_0/F$ } iff\\
 - $A/A_0$ is finite\\
 - $A_0/F$ is periodic and $\p$ is multiplication by $\alpha$, on $A_0/F$.\\
 - $F$ is free abelian and $\p=0$ on $F$.\\
Clearly, we are considering precisely those \emph{inertial endomorphisms $\p$ which are FM on some $A/F$ and have periodic image} (which is not necessarily finite as in Proposition \ref{G}). Observe that such an $\alpha$ is not uniquely determined.

\begin{Prop}\label{RI}  Let $A$ be an abelian group with finite torsion-free rank. Then
the set  $UI(A)$  of uniform inertial endomorphisms of $A$ is a subring of $E(A)$ containing $F(A)$, where:\\
i)\ \ ${UI(A)}/{ F(A)}$ is isomorphic to a subring of ${\cal H}(A)$,\\
ii)\  if $A$ is periodic, $UI(A)=FM(A)=M(A) +  F(A)$ and $UI(A)/F(A)\iso {\cal H}(A)$,\\
iii) there exists a periodic quotient $\bar A$ of $A$ such that
${UI(\bar A)}/{ F(\bar A)}\iso {\cal H}(A)$.
\end{Prop}

\pf  For the definition of ${\cal H}(A)$ see the introduction. Clearly $F(A)\subseteq UI(A)\subseteq IE(A)$.
Consider then  the relation $\cal R$ between $\p\in UI(A)$ and $\alpha\in\J$ defined by the following:
$\p$ is uniform and is represented by $\alpha$ on some $A_{0}/F$.

We claim that \emph{$\cal R$ is
compatible with ring operations}. In fact if
 $\p_i$ is represented by $\alpha_i\in\J$ on $A_i/F_{i}$ ($i=1,2$), then $\p:=\p_1-\p_2$ acts as  multiplication by $\alpha:=\alpha_1-\alpha_2$ on $A_3/F_3$ with
$A_3=(A_1\cap A_2)$ and $F_3:=(F_1+F_2)\cap A_3$.
Hence $\p$ is  uniform and represented by  $\alpha$ on $A_{0}/F$, where $F:=F_1\cap F_2$ and $A_{0}/F:=ker(\p-\alpha)_{|A/F}$ (note that
 $F$ is free, $A/F$ is periodic and $\p(F)=0$, as all $A_i$'s and $F_i$'s are commensurable, resp.). For the multiplicative ring operation
argue the same way. In particular \emph{we have that $UI(A)=dom\ \cal R$ is a subring}.
Analogously, $cod\ \cal R$ is a subring of  $\J$. In the case $A$ is periodic, it is plain that $cod\ {\cal R}={\cal J}$.

On one hand we have $\{\p\in UI(A)\ |\ \p{\cal R}0 \}=F(A)$;
in fact  it is plain
 that if $\p\in F(A)$, then $\p{\cal R}0$. Similarly, if  $\p{\cal R}0$,
then $\p=0$ on some $A_{0}/F$ and on $A/T$ as well. Therefore $\p=0$ on $A_{0}$ (as $F\cap T=0$) and $\p\in F(A)$.

Fix now $F$, related sequences $({e_p})$,
$({\e_p})$ and note that  the following ideal  does not depend on the choice of $F$
 $$I:=\bigoplus_p (p^{\e_p}) + \prod_p (p^{e_p})$$
Let us show that $I=\{\alpha\in\J\ \ |\ 0{\cal R}\alpha \}$. If $0{\cal R}\alpha$ then $p^{\e_p}$ divides
$\alpha_p$ for all $p$ and even $p^{e_p}$ divides
$\alpha_p$ for all but finitely many $p$. Thus $\alpha\in I$.
Conversely, let $\alpha=(\alpha_p)\in \prod_p (p^{e_p})$. Then $0{\cal R} \alpha$ since $\alpha=0$ on  $A/F$
(as the $p$-component of $A/F$ is bounded by $p^{e_p}$). Similarly if $\alpha\in \bigoplus_p (p^{\e_p})$, let
$A_0/F$ be generated by  the $p^{\e_{p}}$ socle of $A/F$ for all $p$ such that $\alpha_p\not= 0$ and the whole $p$-component of $A/F$ for the remaining $p$. Thus  $0{\cal R} \alpha$ as $\alpha=0$ on  $A_{0}/F$ and $A_{0}$ has finite index in $A$.

From what we proved above, it follows clearly that \emph{$\cal R$ induces an isomorphism as in the statement} $(i)$. For statement $(ii)$ see Proposition \ref{ENDO_mod_F}.

For statement $(iii)$, for each $p$ choose an $L_p\ge F$ such that factor $A/L_p$  is either of type $\Z(p^\infty)$ if $p^{\e_p}=\infty$ or an infinite homogeneus $p$-group with bound  $p^{\e_p}$ plus a group of type $\Z(p^{e_p})$, otherwise. Then set $\bar A:=A/\cap_p L_p$ and  get ${\cal H}(\bar A)={\cal H}(A)$. Then $(iii)$ follows from $(ii)$.\qed

We state a Corollary to our Main result which gives \emph{instances of uniform inertial endomorphisms, which are $0$ on both $A/T$ and $T$ but are not finitary}.

\begin{Cor}\label{G} There exists an abelian group $A$ with $r_0(A)=1$ and $p$-components with order $p$ such that the ideal of inertial endomorphisms annihilating $A/T$
 has shape $UI(A)=X+F(A)$, where $X$ is the ideal of (all) endomorphisms annihilating both $A/T$ and $T$.

  Moreover, concerning additive groups, we have $X \iso \prod_p\Z(p)$ and $UI(A)/F(A)$ is an uncountable torsion-free divisible abelian group.
\end{Cor}

\pf Consider the group $A$ as in Proposition A of \cite{DR2}, which has no-critical primes.
On one hand $X \subseteq UI(A)$, by that Proposition. On the other hand,
by our Main Result, $UI(A)$ is the ideal of inertial endomorphisms annihilating $A/T$, as stated.
Arguing now on $T$ as in Remark 4 of \cite{DR3}, we see that any $\p$ in $UI(A)$ acts as a finitary automorphism $\p_{|T}$ on $T$.
Thus there is a finite set $\pi$ of primes such that $\p_{|T}=m\oplus 0$ on $T=T_\pi\oplus T_{\pi'}$. Then there is $C$ such that
$A=T_\pi\oplus C$ and one can consider $\p_1:=m\oplus 0\in F(A)$ where $\p-\p_1$ annihilates both  $T$ and $A/T$. Thus $UI(A)=X+F(A)$. Notice that $X\iso Hom(A/T,T)\iso  \prod_p\Z(p)$ and $X\cap F(A)=T(X)$.
\qed


\section{Proof of Main Result}\label{Proof}
Statement $(a)$ follows immediately from the Characterization Theorem. Let then $A$ have FTFR. Fix sequences $\pi_i=\{p_1,...,p_i\}$, $(B_i)$, $(C_i)$ as in Lemma \ref{mini-moltiplicazioni}. Let $N$ be the ring  of mini-multiplications w.r.t. fixed sequences.

 For each $\p\in IE(A)$, in the notation of the Characterization Theorem, $\p=m/n$ on $A/T(A)$ while $A_\pi$ is bounded where $\pi:=\pi(n)$. Moreover
 $A=A_\pi\oplus_\pi C_*$ where $C_*$ is $\Q^{\pi}$-module, as $T(C_*)$ is $\pi$-divisible being $\pi'$-group and $C_*/T(C_*)\iso A/T(A)$ is $\pi$-divisible as well. Thus there exists the semi-multiplication $0\oplus_\pi m/n$. We may consider also $\p_1:=\p-(0\oplus_\pi m/n)$, which is $0$ on $A/T$. We reduced to show that if $\p_1=0$ on $A/T$ then $\p_1\in UI(A)+N$.

Apply now part $(b)$ of the Characterization Theorem to $\p_1$ and fix $$A_0=B\oplus D\oplus C$$ and $V\le C$ as in that statement.
Then for some sufficiently large $j$ we have $\pi(B)\cap\pi_c(A)\subseteq\pi_{j}\subseteq \pi_c(A)$. Fix $j$ as well. As in Lemma \ref{mini-moltiplicazioni} we have $$A=B_j\oplus_{\pi_j}  C_{j}.$$
Let now, for each $i$, \ $s_i:=0$
if either the $p_i$-component of $(D+C)/V$ if finite or that of $B$ is finite. Otherwise, let $s_i:=m_i-\alpha_i$ where $\alpha_i\in Q_{p_i}^*$ represents $\p_1$ on the $p_i$-component of $(D+C)/V$ and  $m_i\in\Z$ represents $\p_1$ on the $p_i$-component of $B$. Let $s\in \Z$ such that the mini-multiplication $\sigma:=s\oplus_{\pi_j}0$ on $B_j\oplus_{\pi_j}  C_{j}$ acts as multiplication by $s_{i}$ on the $p_i$-component of $B_{j}$.

We claim that $\p_2:=\p_1-\sigma\in UI(A)$, where it is clear that $\p_2=0$ on $V_2:=V\cap C_j$. Then let us show that\emph{ $\p_2$ is FM on $A/V_2$}. Since $\p_2$ is inertial, it is multiplication on all but finitely many primary components of $A/V_2$. Thus
it is enough to show that $\p_2$ is FM on each primary component of $A/V_2$. Fix a prime $p$. If $p\not\in\pi_c(A)$ the statement is plain. If
 $p\not\in\pi(B)$, the statement holds as $\p_2$ acts as $\p_1$ on the $p$-component of $A/V_2$. Finally, if $p\in\pi_j$, consider the box decompositions:
 \\ - $A_0=B'\oplus C'$, where $B'$ is the $p$-component of $B$ and $C'$ is $D+C$ plus the $p'$-component of $B$;
 \\ - $A=B''\oplus C''$, where $B''$ is the $p$-component of $B_j$ and $C''$ is $C_j$ plus the $p'$ component of $B_j$

 Note that \emph{$B'$ and $B''$ (resp. $C'$ and $C''$) are commensurable}. This follows from the fact that $B'/(B'\cap B'')$ on one hand is a bounded $p$-group (as a factor of $B'$) and
on the other hand it has finite rank (as isomorphic to a subgroup of $C''$). For the commensurability of $C'$ and $C''$ argue the same way.

Then consider $A':=(B'\cap B'')\oplus(C'\cap C'')$. Now it is enough to note that $\p_2$ is multiplication on the $p$-component of $A'/V'$, where $V':=V\cap A'$, to see that $\p_2\in UI(A)$, as wished.

Thus we have proved that $$IE(A)=SM(A)+UI(A)+N$$
and that $UI(A)+N$ \emph{is the ideal of inertial endomorphisms of $A$ with periodic image}. Therefore the map $\p\mapsto m/n$ with $\p=m/n$ on $A/T$ gives a monomorphism with kernel $UI(A)+N$, while the isomorphism concerning
${UI(A)}/{F(A)}$ is stated in Proposition \ref{RI}.

Finally, if $A$ is periodic the statement is clear once one applies again Proposition \ref{RI}.
\qed


{

}


\end{document}